\documentclass{amsart}

\usepackage[mathcal]{euscript}
\usepackage{mathrsfs}
\usepackage{amsfonts,amssymb,amsbsy,amsmath}
\usepackage{graphicx,color}
\usepackage{epsfig}
\usepackage[all,cmtip]{xy}
\usepackage{verbatim}

\newcommand{\R}{\mathbb{R}}

\newcommand{\C}{\mathbb{C}}
\newcommand{\Z}{\mathbb{Z}}
\newcommand{\id}{\mbox{id}}

\newcommand{\ba}{\begin{array}}
\newcommand{\ea}{\end{array}}

\newcommand{\Lee}{\mathcal{L}}

\newcommand{\co}{conj}

\theoremstyle{plain}
\newtheorem{theorem}{Theorem}
\newtheorem{lemma}[theorem]{Lemma}
\newtheorem{prop}[theorem]{Proposition}
\newtheorem{coro}[theorem]{Corollary}

\theoremstyle{remark}
\newtheorem{remark}{Remark}

\theoremstyle{definition}
\newtheorem*{defi}{Definition}

\begin{document}

\title{There is a unique real tight contact 3-ball}

\author{Fer\.{\i}t \"{O}zt\"{u}rk}
\address{Bo\u{g}az\.{\i}\c{c}\.{\i} \"{U}n\.{\i}vers\.{\i}tes\.{\i}, Department of Mathematics, TR-34342
  Bebek, \.Istanbul, Turkey}
\email{ferit.ozturk@boun.edu.tr}

\author{Nerm\.{\i}n Salep{c}\.{\i}} 
\address{Mathematisches Forschungsinstitut Oberwolfach,
Schwarzwaldstr. 9-11 (Lorenzenhof),
77709 Oberwolfach-Walke,
Germany}
\email{salepci@mfo.de}

\begin{abstract}
We prove that there is a unique real tight contact structure on the  3-ball with convex boundary 
up to   isotopy through real tight contact structures. We also give a partial
classification of the real tight solid tori with
the real structure being antipodal map along longitudinal and the identity along meridional direction.
For the proofs, we use the real versions of contact neighborhood theorems
and the invariant convex surface theory in real contact manifolds.
\end{abstract}

\maketitle

\section{Introduction}

Motivated by the complex conjugation on complex manifolds, a {\it real structure}
 on an oriented $2n$-manifold $X$ with boundary (possibly empty) is defined as
an involution $c_{X}$ on $X$ which is orientation preserving 
if $n$ is even and orientation reversing if $n$ is odd  and which has the fixed point set of half dimension, 
if  it is not empty.  We extend the definition of a real structure to smooth oriented $2n-1$ dimensional manifolds. 
A {\it real structure} on an oriented  $(2n-1)$-manifold is
an involution which is orientation preserving 
if $n$ is even and orientation reversing if $n$ is odd
and the fixed point set of which is  of  dimension $n-1$, if  it is not empty.
In particular, if $M$ is the oriented boundary of an oriented  $4$-manifold $X$, it can
be shown that a real 
structure on $X$ restricts to a real structure on $M$. For a manifold $W$
with arbitrary dimension, the pair $(W,c_W)$
is called a {\em real (or $c_W$-real) manifold}. The set  fix$(c_W)$ of fixed points of the real structure 
is called the {\em real part} of the real manifold. 

If there is a symplectic (or  a contact) structure on a real manifold
then one can talk about the compatibility of the structure with the real structure.
As the simplest example, consider $B^4$ in $\C^2$ with the real structure given by
the complex conjugation $\co$ and the symplectic structure given by the
standard symplectic form $\omega_{std}$. Observe that $\co$ is  orientation preserving and
$\co^*\omega_{std}=-\omega_{std}$. On $S^3=\partial B^4$, 
the map $c_{S^3}=\co|_{S^3}$ is an orientation preserving involution. It is straightforward to verify
that the induced contact form $\alpha_{std}$ on $S^3$ satisfies $c^*_{S^3}\alpha_{std}=-\alpha_{std}$ too.
This toy example exhibits a general fact. Assume that $(M,\alpha)$ is a contact $(2n-1)$-manifold
which is strongly symplectically fillable by $(X,\omega)$ and that there is a real structure
$c_{X}$ on $X$ preserving $M$ and satisfying ${c^*_{X}}\omega=-\omega$. Then ${c^*_{X}|_{M}}\alpha=-\alpha$.
Thus we are motivated to give the following
\begin{defi}
\begin{enumerate} 
\item (due to C.~Viterbo  \cite{vi}) Let $(X,\omega)$ be a symplectic manifold of dimension $2n$ and let $c_{X}$ be a
real structure on $X$.
If $c^*_{X}\omega=-\omega$, the form $\omega$ is called a 
{\em real  (or $c_X$-real) symplectic form} with respect to $c_{X}$ and
$c_{X}$ is called an {\em anti-symplectic involution}. The triple $(X,\omega, c_{X})$ is called a
{\em real symplectic manifold}.
\item Let $(M,\xi)$ be a contact manifold of dimension $2n-1$ with contact distribution $\xi$. Let $c_{M}$ be 
a real structure on $M$.
If $({c_{M}})_*\xi=-\xi$, the contact structure  $\xi$ is called a {\em real  (or $c_M$-real) contact structure} 
with respect to $c_{M}$.  The triple $(M,\xi,c_{M})$ is called a {\em real contact manifold}.
\end{enumerate}
\end{defi}

In this note, we always assume
that $\xi$ is cooriented so that it is globally defined by the kernel of a contact form
$\alpha$ and hence the real condition reads ${c^*_{M}}\alpha=-\alpha$ as well. 
Moreover,  if  $M$ is oriented we assume that  $\xi$  is positive, viz. the orientation defined by $\alpha\wedge d\alpha$ agrees with  the given orientation of $M$. \\

Real algebraic varieties with the symplectic form induced from the standard one on
the ambient projective space are natural examples for real symplectic manifolds
while links of real algebraic singularities constitute natural examples for real contact manifolds.

In dimension 3, two standard examples for real contact manifolds are \\
(i) $(B^3,\xi_{std},c_{std})$ with $\xi_{std}$ 
given by $\alpha_{std}=dz+xdy-ydx$ and $c_{std}$ rotation by $\pi$ around $y$-axis;\\
(ii) $(T=S^1\times D^2,\eta_{std},c_2)$ with $c_2:(\theta,(u,v))\mapsto (\theta,(-u,-v))$
and $\eta_{std}$ 
given by the 1-form $\cos \theta du - \sin \theta dv$ ($u$, $v$ rectangular coordinates).\\

The real $c_{std}$ in the first example is known to be the unique real structure on $B^3$
up to isotopy through real structures (by P.A.~Smith's work in 1930's).
Whereas on a solid torus, there are four real structures
up to isotopy through real structures  \cite{ha}. Let us consider a solid torus $S^1\times D^2$.
We choose  an oriented  identification of $T^2$ with $\R^2_{(x,y)}/ \Z^2$, where $x$ 
direction corresponds the meridional direction of $T^2$ and fix the coordinates of 
$S^1\times D^2$  as  $(y, x, t)$ where $t$ is the radius direction of $D^2$ 
(in accordance with the convention in \cite{ho}).
In these coordinates, the four possible real structures on the solid torus are:

\noindent (1) $c_1:(y,(x,t))\mapsto (-y,(-x,t))$;\\
(2) $c_2:(y,(x,t))\mapsto (y,(x+\frac{1}{2},t))$;\\
(3) $c_3:(y,(x,t))\mapsto (y+\frac{1}{2},(x,t))$;\\
(4) $c_4:(y,(x,t))\mapsto (y+\frac{1}{2},(x+\frac{1}{2},t))$.

Moreover any orientation
preserving involution on $T^2$ can be extended to an involution on $S^1\times D^2$. Such an extension is
unique up to isotopy and fixes a core of the solid torus setwise \cite{ha}. We will refer the real structures
$c_1,c_2,c_3,c_4$ frequently in the sequel. Likewise there are three involutions on the core circle 
up to isotopy through involutions: reflection, identity and rotation by $\pi$ (antipodal map). 
We also denote them by $c_1$, $c_2$ (or id) and $c_3$ respectively. Note that $c_4|_{\mbox{core}}=c_3$ too.

Before we state the main result of the present work, recall that
up to isotopy through tight contact structures, there is a unique tight structure on
$B^3$ and $\xi_{std}$ is a representative \cite{el}. 
With respect to $\xi_{{std}}$, $\partial B^3$ 
is convex and a dividing curve on it is the equator $\mathcal{E}=\{x^2+y^2=1\}$.
Also we remarked above that $c_{std}$ is the unique real structure on $B^3$.
However, it is not at all obvious that any real tight structure on $B^3$ is
isotopic to $\xi_{{std}}$ through real tight structures.
We prove the uniqueness:

\begin{theorem}\label{ana}
Up to isotopy  through real contact structures making $\partial B^3$ convex, 
there is a unique real tight contact structure  on $B^3$.
\end{theorem} 

To  prove the theorem, we divide the 3-ball into real tight pieces. All are {\em standard} 
except one which is a $c_3$-real tight solid torus. Then we use our partial classification of $c_3$-real tight
structures on the solid torus to conclude.

In Section~2, we prove real contact neighborhood theorems for real isotropic and contact
submanifolds in real tight manifolds of arbitrary odd dimension. To do that we follow
the classical proofs. In Section~3, we prove fundamental theorems of 
invariant convex surface theory in real tight 3-manifolds following \cite{gi}. 
We devote Section~4 to the study of $c_3$-real tight solid tori
where we present a partial result towards the classification of $c_3$-real tight solid tori up to isotopy through real contact structures.
The proof of Theorem~\ref{ana} is also included in this section. The main approach in this section 
is convex surface theory involving dividing sets, bypass attachments and edge-rounding (see \cite{gi}, \cite{ho}).

In what follows,  when we say that a map $f:(M,c)\rightarrow (M',c')$ is {\em equivariant}
we mean  $f\circ c=c'\circ f$. In particular  $f$ is {\em $c$-equivariant} means  $f\circ c=c\circ f$.  Besides, when objects are sent to themselves under a real structure $c$, we call them 
{\em invariant}  (or {\em  $c$-invariant}) or sometimes {\em  symmetric}. 
We say {\em antiequivariant} or {\em antisymmetric}  (or sometimes  {\em $c$-antisymmetric})  if they are sent to {\em minus} themselves. 
For example a submanifold might be invariant (symmetric)
as a set while with an orientation it might be antisymmetric.

\section{Real contact neighborhood theorems}

In this section, we give several real contact neighborhood theorems. 
Although we will not use all in the present work, 
we present them to refer in our future work. We first observe that
{\it Gray stability} is still valid in real setting. In its proof
one can still employ Moser trick, which is valid in real setting too. Indeed, 
to find an equivariant isotopy $\psi_t$ of  contact forms, we consider it as the flow of a time-dependent
antiequivariant vector field $X_t$. Since the differential equation for $X_t$
is antiequivariant, the solution for $X_t$ must be symmetric (see e.g. the proof of \cite[Theorem~2.2.2]{ge}).
It is easy to check that an appropriate version of 
Equivariant Darboux-Weinstein Theorem (\cite[Theorem~22.1]{gs})
for symplectic manifolds
is valid in the real contact setting. To do that, we go through the proofs of contact
neighborhood theorems, minutely presented in \cite{ge}, and make sure that the real
versions of those theorems are valid. Below, we do not repeat the proofs; instead we point out 
where the proofs in op.cit. have to be altered to apply to our setting. Hence
the following discussion is not self-contained so it must be read in parallel with the 
related proofs in \cite{ge}.

In the sequel we will use the notion of a real symplectic bundle the definition of which is given here:
\begin{defi} 
Let $(E,B,\pi,\omega)$ be a $2k$-dimensional  symplectic vector bundle over a manifold $B$ and the bundle map $\pi$.
Assume there is an involution $c_E$ on $E$ which is a  bundle isomorphism covering the identity and which gives a real structure on each fiber.
The triple $(E,\omega,c_E)$ is called a {\em real symplectic bundle}.

\end{defi}

\begin{theorem} (Real Contact Neighborhood Theorem for Real Parts)
\label{LegKomsu}
Let $(M_i,\xi=\ker \alpha_i,c_i),\; (i=0,1)$, be two real contact manifolds, and $L_i$ closed isotropic
submanifold of $M_i$ lying in the real part of $M_i$. Suppose that an equivariant real symplectic bundle 
isomorphism 
$$\Phi:((TL_0)^{\perp_{d\alpha_0}}/TL_0,d\alpha_0,Tc_0)\rightarrow ((TL_1)^{\perp_{d\alpha_1}}/TL_1,d\alpha_1,Tc_1)$$
is given and that $\Phi$ covers a diffeomorphism  $\phi:L_0\rightarrow L_1$.
Then there exist neighborhoods $\nu(L_i)$ of $L_i$ and a diffeomorphism 
$F:\nu(L_0)\rightarrow \nu(L_1)$ such that:\\
(i) $F|_{L_0}=\phi$;\\
(ii) $F^*\alpha_1 =\alpha_0$;\\
(iii) $F$ is equivariant.
\end{theorem}

\noindent {\it Proof:} We follow the proof of \cite[Theorem~2.5.8]{ge}.  
Recall that (dropping the index $i$) there is a decomposition
$$ T_L M = \left<R_{\alpha}\right> \oplus TL \oplus (\xi|_L)/(TL)^{\perp_{d\alpha}} \oplus (TL)^{\perp_{d\alpha}}/TL$$
where $R_{\alpha}$ is the Reeb vector field of $\alpha$ and $\perp_{d\alpha}$ denotes the orthogonal complement
with respect to $d\alpha$. We first observe that $\left(TL\oplus (\xi|_L)/(TL)^{\perp_{d\alpha}},d\alpha, Tc\right)$  is a real symplectic bundle.
Also consider the real symplectic bundle  $(TL\oplus T^*L,\Omega_{L}, s)$

where $s=\id \oplus -\id$ and $\Omega_{L}$ is defined by:
$$(\Omega_{L})_p (X+\eta,X'+\eta') = \eta(X')-\eta'(X); \;\; X,X'\in T_p L, \eta,\eta'\in T^*_p L.$$
It is straightforward to confirm that  $s(\Omega_{L})=-\Omega_{L}$.

Now recall (see \cite[Lemma~2.5.7]{ge}) that there is a symplectic bundle isomorphism
\begin{equation} 
\ba{rl}
\id_{TL}\oplus \Psi: \left(TL\oplus (\xi|_L)/(TL)^{\perp_{d\alpha}}, Tc, d\alpha\right) 
& \rightarrow (TL\oplus T^*L,\Omega_{L}, s),\\ 
(X,Y) & \mapsto (X,(\iota_{Y}d\alpha)|_{TL})
\ea
\label{denk}
\end{equation}

Thus having symplectic bundle isomorphisms onto $TL_i\oplus T^*L_i$ for each $i$, and using the
symplectic bundle isomorphism $\Phi$ for the remaining symplectic part of $T_{L_i} M_i$,
one can build a symplectic bundle isomorphism from $T_{L_0} M_0$ to $T_{L_1} M_1$.
This finishes the linear algebra part of the proof in the lack of a real structure.

In the real setting, we first check the validity of \cite[Lemma~2.5.7]{ge}. 

\begin{lemma} 
The bundle map in (1) is  an equivariant isomorphism of real symplectic bundles.
\end{lemma}

\noindent {\it Proof:} We just need  to check that the map $\id_{TL}\oplus \Psi$ is equivariant. In fact,
$$\ba{rcl} 
(\id_{TL}\oplus \Psi)_p\circ (Tc)_{p} (X,Y) & = & (\id_{TL}\oplus \Psi)_p (X,(Tc)_{p}(Y))\\
& = & \left(X,(\iota_{(Tc)_{p}(Y)}d\alpha)|_{T_p L}\right) \\
& = & \left(X,d\alpha|_{T_p L}\left((Tc)_{p}(Y),(Tc)_{p}(^{\centerdot})\right)\right) \\
& = & \left(X, (\iota_{Y} c^* d\alpha)|_{T_p L}\right) \\
& = & \left(X, -\iota_{Y} d\alpha|_{TL}\right)_p \\
& = & \left(s\circ \id_{TL}\oplus \Psi\right)_p (X,Y)
\ea
$$
For the third equation, we use the facts that $d\alpha$ is restricted on $TL$ and $Tc$  restricted to $TL$ 
is the identity map, $L$ 
lying in the real part. The fifth equation follows from the fact that $d\alpha$ is real with respect to $c$.
\hfill  $\Box$ \\

Since $\phi$ is equivariant, we have already  showed that the map
$$ \ba{rl}
T\phi \oplus (\Psi_1^{-1}\circ (\phi^*)^{-1}\circ \Psi_0): &
\left(TL_0\oplus (\xi|_{L_0})/(TL_0)^{\perp_{d\alpha_0}}, d\alpha_0, Tc_0\right) \\
 & \rightarrow \left(TL_1\oplus (\xi|_{L_1})/(TL_1)^{\perp_{d\alpha_1}}, d\alpha_1, Tc_1\right) 
\ea
$$
is an equivariant isomorphism of real symplectic vector bundles.

Let $\Phi_R: \left< R_{\alpha_0}\right> \rightarrow \left< R_{\alpha_1}\right>$
be  the bundle map taking $R_{\alpha_0}(p)$ to $R_{\alpha_1}(\phi(p))$. It is 
equivariant, i.e. $\Phi_R\circ Tc_0 = Tc_1 \circ \Phi_R$; in fact $Tc_i(R_i)=-R_i$,
since $\alpha_i$ is $c_i$-real.
Then the bundle map
$$ 
\ba{rl}
\tilde{\Phi} = 
\Phi_R\; \oplus & \!\!\!(\Psi_1^{-1}\circ (\phi^*)^{-1}\circ \Psi_0)\; \oplus \; \Psi:  \\
&\left<R_{\alpha_0}\right> \oplus TL_0 \oplus (\xi|_{L_0})/(TL_0)^{\perp_{d\alpha_0}} \oplus (TL_0)^{\perp_{d\alpha_0}}/TL_0 \\
&\longrightarrow 
\left<R_{\alpha_1}\right> \oplus TL_1 \oplus (\xi|_{L_1})/(TL_1)^{\perp_{d\alpha_1}} \oplus (TL_1)^{\perp_{d\alpha_1}}/TL_1
\ea
$$
is an equivariant bundle isomorphism of real bundles, preserving the structures induced from the contact structure
on each subspace. Hence we are done with the linear algebraic part of the proof.

Now we extend this
real symplectic vector bundle isomorphism to neighborhoods of $L_0$ and $L_1$.
As in the proof of \cite[Theorem~2.5.8]{ge}, we choose similar tubular maps $\tau_i:(NL_i,Tc)\rightarrow (M_i,c)$
requiring one more condition that $\tau_i$ is equivariant. This can be achieved by choosing an equivariant
Riemannian metric in a neighborhood of $L_i$ in $M_i$. Hence we obtain
a diffeomorphism $\tau_1 \circ \tilde{\Phi} \circ \tau_0^{-1}:\nu(L_0)\rightarrow \nu(L_1)$
which satisfies the conditions {\em (i), (iii)} of the theorem on $L_i$'s. Finally, the pull-back of
$\alpha_1$ in $\nu(L_0)$ can be equivariantly isotoped to $\alpha_0$ thanks to the
real version of Gray stability. Composing the time-1 diffeomorphism with the diffeomorphism above,
condition {\em (ii)} is satisfied too.\hfill $\Box$ \\

In particular, since the 1-dimensional real part in any contact real 3-manifold is Legendrian, we obtain

\begin{coro} (Real Contact Neighborhood Theorem for Real Knots)
\label{gerceldugum}
A real knot in an arbitrary real contact 3-manifold
has  a neighborhood equivariantly contactomorphic to the real contact 3-manifold
$(S^1\times D^2,\eta_{std},c_2)$, defined in the Introduction, 
with $L$ having the image  $S^1\times\{0\}$.
\end{coro}

We now state a similar theorem for neighborhoods of real knots up to {\em equivariant contact isotopy}.
Its proof lies in the proof of Theorem~\ref{LegKomsu} and essentially
the same as the proof for the standard neighborhood theorem for
Legendrian knots up to contact isotopy. The idea is to classify the
real plane bundles over the knot as in the proof of Theorem~\ref{LegKomsu} and then to use the 
exponential map corresponding to an equivariant Riemannian metric.

We denote by $\mbox{tw}(L,s)$ the twisting number of $\xi$ along $L$
with respect to a fixed trivialization $s$ of the normal bundle of the Legendrian knot
$L$ in $M$. Similarly we can define the twisting number of a Legendrian arc relative end points. 
Note that in that case  $\mbox{tw}(L,s)$ would be a real number. If $s$ is understood
or its omission from the notation does not affect the discussion, we write
simply $\mbox{tw}(L)$.

\begin{coro} (Real Contact Neighborhoods of Real Knots up to Equivariant Contact Isotopy)
\label{gercelyay}
Let $L$ be a real knot  in a real contact 3-manifold $(M,\xi,c)$. Given an identification
with $S^1\times D^2$ and the longitudinal framing along $L$,
an invariant neighborhood of $L$ is  isotopic through real contact structures 
to the real contact 3-manifold
$$(S^1\times D^2,\cos(\mbox{tw}(L) \theta) du + \sin( \mbox{tw}(L) \theta)  dv, c_2)$$ 
where $(\theta,(u,v))\in S^1\times D^2$. If $L$ is a real arc,
then  an invariant neighborhood of $L$ is isotopic relative end points 
through real contact structures to

$$([0,1]\times D^2,\cos(2\pi\mbox{tw}(L)\theta) du + \sin(2\pi\mbox{tw}(L)\theta)  dv, c_2).$$
\end{coro}

\noindent {\it Proof:}
The proof of Theorem~\ref{LegKomsu} shows that isotopy classification of the contact structures 
near $L$ turns out to be topological so  
the first part of the corollary is obvious. As for the second part, we need to explain how to construct a smooth neighborhood of  the arc $L$.

Consider an auxiliary Legendrian arc $L'$ with end points $p_{0}$ and $ p_{1}$ (such an $L'$ exists). 
The union $L\cup L'$  is a Legendrian knot  whose contact neighborhood is isotopic 
to ($S^1\times D^2$,  $ \cos( \mbox{tw}(L\cup L') \theta) du + \sin( \mbox{tw}(L\cup L') \theta)dv , c_{2})$ 
(Corollary~\ref{gerceldugum}). We can assume that by a rescaling of $L\cup L'$ 
the contact planes  $\xi(p_{0})$ and $\xi(p_{1})$ stay fixed.  
Consider the (non-smooth) neighborhood  $\nu(L)$ of $L$ 
given by $[p_{0}, p_{1}]\times D^2$.  We need to  smooth out the corners of $\nu(L)$ 
in such a way that the boundary of $\nu(L)$ becomes convex. 

We can work in the following model for  $\partial\nu(L)$ around  $\{p_{1}\}\times D^2$.   
In $(R^3, \xi_{std}, c_{std})$, consider the surface $S$ which is the
union of the cylinder $\{x^2+z^2=r^2,$  $y\in [1-\epsilon, 1]\}$ 
and  the disk $\{x^2+z^2\leq r^2,  y=1\}$.   It is straightforward to check that 
$X= \partial_{y} +(x+kz)\partial_{z}+kx\partial_{x}$ ($k\in {\R}$) is  an 
invariant contact  vector field  and 
transverse to $S$ for sufficiently large  $k\in {\R}^+$. Moreover, $X$ is also transverse to the 
smooth surface obtained by rounding the edge of  $S$ with the  {\em quarter}  torus obtained by rotating 
about the $y$-axis the quarter circle
$\{(y-1+\delta)^2+ (z-r+\delta)^2=\delta^2, y\geq 1-\delta, z\geq r- \delta\}$.  Hence the
smoothened surface is convex.
Note also that the construction  is invariant with respect to the real structure $c_{std}$.
\hfill  $\Box$ \\

One can prove a theorem similar to Theorem~\ref{LegKomsu} for invariant contact submanifolds. We will
use an equivariant version of the set-up for \cite[Theorem~2.5.15]{ge}.
Let $(M,\xi=\ker \alpha,c)$ be a real contact manifold and $(K,\xi|_{K})$ a closed invariant
contact submanifold of $M$. 
Then we have $$ TM|_K = TK \oplus NK = TK \oplus (TK)^{\perp_{d\alpha}}/TK$$
as smooth vector bundles.

\begin{theorem} (Real Contact Neighborhood Theorem for Invariant Contact Submanifolds)
Let $(M_i,\xi_i,c_i),\; (i=0,1)$ be two real contact manifolds, and $(K_i,\xi_i|_{K_{i}})$ closed invariant
contact submanifolds of $M_i$. Suppose there is an isomorphism of equivariant conformal symplectic normal
bundles 
$$\Phi:((TK_{0})^{\perp_{d\alpha_{0}}}/T_{0}, Tc_{0}) \to ((TK_{1})^{\perp_{d\alpha_{1}}}/T_{1}, Tc_{1})$$ covering a given
equivariant contactomorphism $\phi:(K_0,\xi_0|_{K_{0}})\rightarrow (K_1,\xi_1|_{K_{1}})$. 
Then there exist neighborhoods $\nu(K_i)$ of $K_i$ and a diffeomorphism 
$F:\nu(K_0)\rightarrow \nu(K_1)$ such that:\\
(i) $F|_{K_0}=\phi$;\\
(ii) $F^*\alpha_1 =\alpha_0$;\\
(iii) $F$ is equivariant;\\
(iv) $TF|_{(TK_{0})^{\perp_{d\alpha_{0}}}/T K_0}$ and $\Phi$ are equivariantly bundle homotopic as conformal symplectic bundle
isomorphisms.
\end{theorem}

\noindent {\it Proof:} As in the proof of Theorem~\ref{LegKomsu} above,
it is enough to build an equivariant bundle map from $TM_0|_{K_0}$ to $TM_1|_{K_1}$
covering $\phi$ and inducing $\Phi$.

For this it suffices for us to secure that the choices made in the proof of \cite[Theorem~2.5.15]{ge}
can still be made in the real setting. The rest of the proof is exactly the same as in the proof
of that theorem. Now, we go into that proof.
We first choose a real contact form $\alpha'_i$ for ${\xi_i}|_{K}$ ($i=0,1$).
Whatever the choice is, we want to scale $\alpha'_i$ to a {\em real} contact
form  $\alpha_i$ such that $\alpha_i$ gives 1 when contracted 
with the Reeb vector field $R'_i$ corresponding to $\alpha'_i$. The scaling function $f_i$ is found to
satisfy $df_i=\iota_{R'_i} d\alpha_i$ on $TM_i|_{K_i}$. Since $R'_i$ and  $\alpha_i$ are antisymmetric,
the function $f_i$ is found to be equivariant.
\hfill  $\Box$ \\

\section{Invariant convex surfaces}

Let $S$ be a closed connected oriented surface in a real 3-manifold $(M,c)$ such that $c(S)=S$.
If $c|_{S}$ is orientation preserving (i.e. $S$ is symmetric), 
then the real points of $S$ (if any) are isolated whereas if
$c|_{S}$ is orientation reversing (i.e. $S$ is  antisymmetric), 
then the set of real points on $S$ is a disjoint union (possibly empty) of circles.

\begin{theorem}
\label{eqVf}
Let $S$ be an oriented, symmetric (antisymmetric) convex surface in $(M,\xi,c)$.
Then there exists a symmetric (respectively antisymmetric) contact vector field for $\xi$ near $S$
and thus an invariant dividing set on $S$.
\end{theorem}

\noindent {\it Proof:} Since $S$ is convex, there is a contact vector field $Y$ over $S$ by definition,
i.e. $ \Lee_{Y} \alpha = \mu \alpha$  for some contact form $\alpha$ defining $\xi$ in a neighborhood $\nu S$
of $S$ and some function $\mu:\nu S\rightarrow \R$.
We claim that  if $c|_S$ preserves (reverses) the orientation of $S$ then the symmetric vector field
$X=Y+c_*Y$  (respectively the antisymmetric vector field $X'=Y-c_{*}Y$) is a contact vector field. 
Below, calculations are done for  $X$, the case of $X'$ is similar.

If  $c|_S$ preserves the orientation then $Y$ and $c_{*}Y$ are on the same side of $S$. 
This assures that $Y+c_{*}Y$ is never zero and is transverse to $S$.
Moreover 
we have:
$$ \ba{rcl} \Lee_{Y+ c_*Y} \alpha& = & \Lee_{Y} \alpha + \Lee_{c_*Y} \alpha \\
&=  & \mu\alpha +c^* \Lee_{Y} c^*\alpha \\
&=  & \mu \alpha - c^*(\mu\alpha) \\
&=  &  (\mu + \mu\circ c)\alpha,
\ea $$
where  $\mu + \mu\circ c:\nu S\rightarrow \R$.

Let $\Gamma\subset S$ denote the dividing set for $X$.  
By definition, $p \in \Gamma$ if and only if $X_{p}$ is tangent to  $\xi_{p}$ or equivalently $\alpha_{p}(X_{p})=0$.
The set $\Gamma$ is invariant  if $\alpha_{c(p)}(X_{c(p)})=0$. Since $c^*\alpha=-\alpha$, 
we have $0=\alpha_{p}(X_{p})=-(c^*\alpha)_{p}(X_{p})=-\alpha_{c(p)}c_{*}(X_{p})$.   
It was shown above that $(c_{*}(X))_{p}= X_{c(p)}$, thus $c(\Gamma)=\Gamma$.
\hfill  $\Box$ \\

\begin{remark}\label{orienondiv}
Recall that the characteristic foliation of a contact structure inherits an orientation from
the orientation of $S$ and the coorientation of the contact structure $\xi$.
Namely, suppose $\Omega$ is a volume form on $S$ and $\alpha$ a 1-form defining $\xi$.
The characteristic foliation of $\xi$ is defined as the integral curves
of the vector field $X$ satisfying $\iota_{X}\Omega=\alpha|_{S}$ \cite[Lemma~2.5.20]{ge}.
This equality determines the orientation of the characteristic foliation.
It follows from the equality $c^{*}\alpha=-\alpha$ that
the real structure reverses (respectively preserves) the orientation of the characteristic foliation
if $c|_S$ is orientation preserving (respectively reversing).
Let us note that the dividing curve is also canonically oriented using the orientation
of the characteristic foliation and of $S$.
Thus, the orientation of the dividing curve is reversed by $c$ regardless of
$c|_S$ being orientation reversing or preserving.  Notice that Theorem~\ref{eqVf} showed 
that the nonoriented dividing set is invariant. \\
\end{remark}

Now we prove an existence theorem for invariant oriented convex surfaces.

\begin{prop}
\label{ahenkcnvx}
Let $S$ be an invariant closed  surface in a real
contact 3-manifold $(M,\xi,c)$. Assume $S$ is oriented and $c|_S$ is orientation preserving.
Then there is a surface $S'$ equivariantly isotopic and $C^{\infty}\!$-close
to $S$ such that $S'$ is convex for $\xi$.
\end{prop}
\noindent {\it Proof:}
Consider the quotient $q:S\rightarrow \bar{S}=S/c$.
Since $c|_S$ is orientation preserving, fix$(c|_S)$ is either empty or
a finite number of points, say $p_1,\ldots,p_{k}$.
%
Let us denote by $X$  the  antisymmetric vector field on $S$
which is the infinitesimal generator
of the restriction of the characteristic foliation $S|_{\xi}$.
Observe that $S|_{\xi}$ does not project
to an orientable foliation on $\bar{S}$; nevertheless, $q_*(S|_{\xi})$ describes a line field
$\mathfrak{F}$ on $\bar{S}$. The fixed points $p_i$ (if they exist)
are regular points of $X$, therefore  each $q(p_i)\in \bar{S}$ is a 
{\it nontrivial singular point}, around which the behavior of $\mathfrak{F}$  is as depicted in
Figure~1.12 of \cite{bn}. Hence our situation satisfies the hypothesis of  \cite[Theorem~2]{bn}.
From its proof  and the openness of the contact condition,
it follows that there is a $c$-antisymmetric Morse-Smale vector field $X_{\xi_1}$ on $S$ arbitrarily close to $X$
which corresponds to a contact structure $\xi_1$,  and which is isotopic to $X$ through 
$c$-antisymmetric vector fields $X_{\xi_t}$, $t\in[0,1], \xi_0=\xi$. Using Gray stability in real setting,
we deduce that there is a surface $S'$ equivariantly isotopic to $S$ such that 
the characteristic foliation $S'|_{\xi}$ is Morse-Smale. 
Since a surface in a contact 3-manifold with  Morse-Smale characteristic vector
field is convex, the proof follows \cite[Proposition~2.16]{gi}.
\hfill $\Box$ \\

In the hypothesis of the theorem above, it is more problematic if we
assume that $c|_S$ is orientation reversing. In that case there are
two main difficulties. First, the quotient $\bar{S}$ may be nonorientable.
The denseness of Morse-Smale vector fields is still an open question
for nonorientable surfaces (see e.g. \cite{pm}). Second, $\bar{S}$ may have boundary. The
boundary curves are in general not Legendrian; nor the field of line 
elements on $\bar{S}$ is in a general position with respect to the boundary.
These make it difficult to conclude the denseness of Morse-Smale vector fields
in our setting.

From the existence of an invariant contact vector field, it follows
immediately
that $S$ has an invariant neighborhood with the contact form
vertically invariant
(see e.g. \cite[Lemma~4.6.19]{ge} and the discussion following it).
Furthermore, we have:

\begin{prop}
\label{GirFlex}
(Real version of Giroux's Flexibility Theorem)
Suppose $S$ is a closed convex surface in $(M,\xi=\ker\alpha,c)$ such that $c(S)=S$;
$S_{\xi}$ the (oriented) characteristic foliation and $\Gamma$ an antisymmetric oriented dividing set.
If $F$ is another oriented  symmetric (or antisymmetric) singular foliation with the same oriented dividing set,
and if the orientations of $F$ and $S_{\xi}$ agree in a neighborhood of $\Gamma$,
then there is an isotopy $\phi_s$ ($s\in [0,1]$) of $S$ satisfying: {\it (i)} $\phi_0=id$,
{\it (ii)} $\phi_s$ is equivariant; {\it (iii)} $\phi_s|_{\Gamma}=id$; {\it (iv)}â $\phi_s(S)$
is convex and {\it (v)} $\xi|_{\phi_1(S)}=F$.
\end{prop}

To prove this version of the theorem, it is enough to observe that every step
of the original proof (see \cite[Theorem~1.2(b)]{gi}) 
can be made $c$-invariant, leading an antisymmetric differential equation
in the standard Moser-type argument. That differential equation has
a unique solution which is necessarily symmetric.

\begin{coro}
\label{kaydir}
The surface $S$ being as above,
let $\xi$ and $\xi'$ be two contact structures with respect to which  $S$ is  convex.
Let $\Gamma$ and $\Gamma'$ be the corresponding  oriented antisymmetric dividing sets on $S$.
If there is an equivariant orientation preserving diffeomorphism 
$f:S\to S$ which is equivariantly isotopic to identity and which sends $\Gamma$ to
$\Gamma'$ preserving the orientation  then $\xi$ and $\xi'$   are 
isotopic through real contact structures
in a sufficiently small neighborhood of $S$.
\end{coro}

\noindent {\it Proof:}
Via the isotopy $\Gamma'$ is made to coincide with $\Gamma$ while $\xi'$ is changed to $\xi''$.
We can assume that in a neighborhood of $\Gamma$ the oriented characteristic foliations for $\xi$ and $\xi'$ coincide, since $f$ preserves the orientation of $S$. Now $\xi$ and $\xi''$  have $\Gamma$ as their common oriented dividing curve. The proof of
Proposition~\ref{GirFlex} shows that $\xi''$ can be isotoped to $\xi$ through real contact structures in a sufficiently small neighborhood of $S$. \hfill $\Box$ \\

\section{Classification of  real tight contact structures on $(S^1\times D^2,c_3)$}

We consider a contact solid torus $S^1\times D^2$ whose boundary $T^2$ is convex.
Let us denote the dividing set of $T^2$ by $\Gamma$ which is determined up to isotopy by the number of
connected components $\#\Gamma$, and their slope  $s$.  We choose coordinates on solid torus as described 
in the Introduction.  
Using the fact that every curve is isotopic to a linear curve on $\R^2/\Z^2$  the slope of the dividing curve is defined as the slope of the corresponding linear curve on $\R^2/ \Z^2$.

In this section, we  prove the following classification theorem
for  real tight contact structures on $S^1\times D^2$ with $\#\Gamma=2$ up to $c_{3}$-equivariant isotopy.
Note that with respect to the above coordinates the real structure $c_{3}$ is
given by  $(y, (x, t))\mapsto ( y+\frac{1}{2}, (x, t)).$

\begin{theorem}
\label{anadolu}
Up to equivariant isotopy through tight structures relative boundary there is a unique real 
tight contact structure on $(S^1\times D^2,c_3)$ having convex boundary with $\#\Gamma=2$ and 
with slope $\pm\frac{1}{2k+1}$ ($k\in\Z, k\geq 0$).
If the slope is $\pm\frac{1}{2k}$ ($k\geq 0$) or $-k$ ($k>1$)
there is no real tight contact structure on $(S^1\times D^2,c_3)$.
\end{theorem}

In the proof of the theorem, we will see that one can peel off  $S^1\times D^2$ by  means of 
$c_3$-invariant double slices and possibly a single basic slice. The existence of 
these  leads an overtwisted structure
except the case the slope equals $\frac{1}{2k+1}, k\in \Z$.

Using the above theorem we can prove the main theorem of this article:\\

\noindent {\it Proof of Theorem~\ref{ana}:}

Without loss of generality we may assume that the real structure on the given  3-ball $B^3$
is $c_{std}$. With this real structure, the arc $L=B^3\cap y$-axis is the real part.

\begin{figure}[h]
\begin{center}
\resizebox{6cm}{!}
{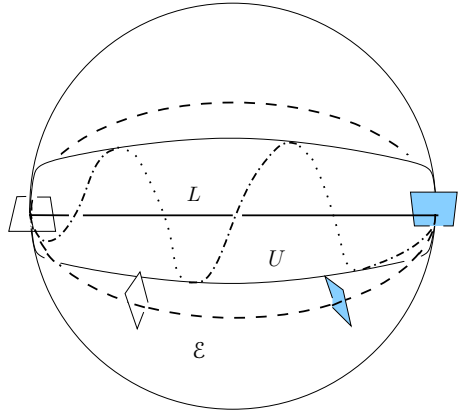}
\caption{The 3-ball  with a neighborhood of its boundary and of the real arc $L$ in standard form.}
\label{topparcala}
\end{center}
\end{figure}

On the convex boundary $S^2$ of $(B^3,c_{std})$ the dividing set $\Gamma$  is  connected
\cite{gi2} and can be chosen symmetric (Theorem~\ref{eqVf}). Any such oriented symmetric
circle on $S^2$ can be equivariantly ambiently isotoped (preserving the orientation) 
to the equator $\mathcal{E}=\{z=0\}$ oriented counterclockwise in $x$-$y$ plane. 
Note that $\mathcal{E}$ is an invariant  dividing set
for  $\xi_{std}$. Thus  we can use Corollary~\ref{kaydir}
to deduce that $\xi$ can be equivariantly isotoped to $\xi'$ in a neighborhood  of
$S^2$ so that in that neighborhood $\xi'$ and $\xi_{std}$ coincide, with the dividing set $\{z=0\}$.
We observe that 
$\mbox{tw}(L,\mathcal{E})$ in $(B^3,\xi')$  is equal to $m+\frac{1}{4}$ for some $m\in \Z$ because of the condition at the end points.
Furthermore, $\xi'$ in  a neighborhood $U$ of $L$ in $B^3$ can be equivariantly
isotoped relative (a neighborhood of) end points so that around $L$,
$\xi'$ has the standard form prescribed  in Corollary~\ref{gercelyay}.
Topological boundary $\partial U$ of $U$ coincides with $S^2$ at small disks $D_{\pm}$ around the real points 
$p_{\pm}=(0,\pm 1,0)$ on $S^2$.
On $D_{\pm}$ the standard form of the contact structure around $L$ and the standard form of
the contact structure around $S^2$ can be made to coincide.
Furthermore on $\partial U$, the dividing set  is connected and {\em turns} $m$ times
around $y$-axis.
In this way, up to equivariant contact isotopy, we have the picture depicted in  Figure~\ref{topparcala}.

Now, $T=B^3 \setminus (int(U)\cup D_{\pm})$
is a solid torus which can be smoothened out about $\partial D_{\pm}$ as follows.
First observe that $\partial D_{\pm}$ is not Legendrian (nor transverse). Therefore
the usual edge-rounding argument (see \cite[Lemma~3.11]{ho}) does not work to make $T$ smooth.
Nevertheless, we can use the model we built in the proof of Corollary~\ref{gercelyay}. 
In $(\R^3, \xi_{std}, c_{std})$ consider the plane $\{y=1$, $x^2+z^2\geq r^2, r\ll 1\}$ (representing a part of $S^2$)
and the cylinder $\{x^2+z^2=r^2, y\in [1-\epsilon, 1], \epsilon\ll r\}$ 
(representing a part of $\partial U\setminus int(D_+)$ near $\partial D_{+}$)
They intersect along the circle $\{x^2+z^2=r^2, y=1\}$. 
We round the edge  by the {\em quarter} torus $\tau$ obtained as the rotation 
about the $y$-axis the quarter circle  
$\{(y-1+\delta)^2+ (z-r-\delta)^2=\delta^2 , y\geq1-\delta , z \leq r+\delta\}$ ($0<\delta<\epsilon$). 
The vector field $X=-s\partial_{y}+ (-sx+kz)\partial_{z} + kx\partial_{x}$, with $s,k\in \R$, is a
symmetric contact vector field with respect to $\xi_{std}$.  
Moreover, $X$ is outwards transverse to the smoothened surface for sufficiently small negative $k$ and $s>0$.

\begin{figure}[h]
\begin{center}
\resizebox{5cm}{!}
{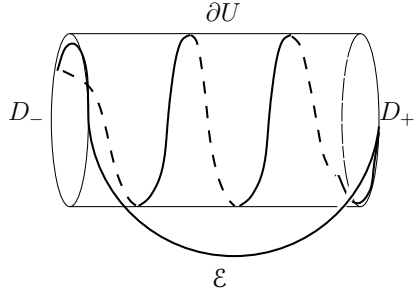}
\caption{The dividing curve
turns $-m$ times around the annulus $\partial U\setminus int(D_{\pm})$. The equator $\mathcal{E}$ is on the outer
annulus $S^2\setminus int(D_{\pm})$}
\label{patates}
\end{center}
\end{figure}

Therefore, $T$ with rounded edges is a real tight  solid torus $(S^1\times D^2,c_3)$ 
with symmetric convex boundary. Furthermore by a straightforward calculation
the dividing set $T$ can be obtained by connecting each component 
of the dividing curve on $\partial U \setminus int(D_+)$ to the  opposite component 
of the dividing curve on $S^2\setminus int(D_+)$ by going right (with $T$ oriented outwards) 
while edge-rounding (same claim for the edge-rounding
near $D_-$). Thus the slope of the dividing curve on $T$ is $-\frac{m+1}{1}$ (see Figure~\ref{patates}).
By Theorem~\ref{anadolu} above, we know that such a solid torus exists  if and only if $m=0$.
Moreover in that case the real tight structure on the solid torus is unique up to isotopy through real contact structures.
 \hfill $\Box$ \\

Theorem~\ref{anadolu} emerges from propositions  and  lemmas presented below. In the following sections
we fix $c_3$ as the real structure and sometimes do not refer it explicitly.
We always assume that $\#\Gamma=2$ on the tori below.

A convex torus in a tight contact manifold has a standard form as follows \cite{ka}.
The dividing set $\Gamma$ has even number of homotopically
essential, parallel curves. In between each pair of dividing curves there is a Legendrian
curve (called  {\em Legendrian divide})  parallel to the dividing curves.
The leaves of the characteristic foliation are transverse to $\Gamma$
and are Legendrian; they define a family of closed Legendrian curves (called
{\em Legendrian rulings}) exhausting the torus.
Thanks to Theorem~\ref{eqVf}, the standard form for tori above can be made symmetric,
so that the dividing set and the set of  Legendrian divides become symmetric. Furthermore,
using the real version of Giroux's Flexibility Theorem (Proposition~\ref{GirFlex})
the  characteristic foliation can be made invariant and
the slope of the Legendrian rulings can have any value except that of Legendrian divides.
We say that such a torus is in {\em symmetric standard from}.

\subsection{Non-existence results for real contact solid tori}

First, let us note that there is no tight contact structure on $S^1\times D^2$ with convex boundary with slope 0 and  $\#\Gamma=2$.
It is simply because each Legendrian divide, being of slope 0, bounds a
meridional disk in $S^1\times D^2$. This disk is an overtwisted disk
since along $L$ the contact planes remain tangent to $T^2$,  so the surface framing and the contact framing are the same,
i.e. the twisting of $L$ is 0 with respect to the meridional disk.

\begin{prop}\label{sonsuz}
There is no real tight contact structure on $(S^1\times D^2,c_3)$
with convex boundary with slope $\infty$ and  $\#\Gamma=2$.
\end{prop}

\noindent {\it Proof:}
Consider a real convex torus in the standard form whose dividing curves have slope $\infty.$
With respect to $c_3$, each component of $\Gamma$ is fixed.
Moreover, this action preserves the orientation of $\Gamma$. But $c_3$ is orientation preserving on $S$
and hence must be reversing on $\Gamma$ (as already discussed above in Remark~\ref{orienondiv}).
 \hfill $\Box$ \\

\begin{remark}\label{equivdt}
A diffeomorphism of $T^2=  S^1\times \partial D^2$ takes meridians to meridians
if and only if it extends to a diffeomorphism of $S^1\times D^2$ (see e.g. \cite[Theorem~2.4]{ha}).
In particular, a Dehn twist along a meridian of $T^2$ extends to a diffeomorphism of  $S^1\times D^2$.  Thus  by
changing  $S^1\times D^2$ using a diffeomorphism, which is a  Dehn twist along a meridian when restricted to $T^2$,
one can change the slope, in particular, make the slope negative. A similar idea applies in the real case.
However, in order not to  spoil the symmetry one should consider a pair of Dehn twists along symmetric meridional curves.
Such a pair of Dehn twists will be called an {\emph{equivariant meridional Dehn twist}}. Note that by means of equivariant Dehn twists 
we can always adjust the slope $s$ so that  $s=-\frac{p}{q}\leq -\frac{1}{2}$ with  $p,q\in \Z^+$, $(p,q)=1$.\\
\end{remark}

\begin{coro}
\label{ciftyok}
There is no real tight contact structure on $(S^1\times D^2,c_3)$
with convex boundary with slope $\frac{1}{2k}, (k\in\Z, k\neq 0)$ and  $\#\Gamma=2$.
\end{coro}

\noindent {\it Proof:}
As it is discussed in Remark~\ref{equivdt},  by sufficiently many equivariant meridional Dehn twists,
we can change the slope to $\infty$ without changing the number of dividing curves.
Hence, the result follows from Proposition~\ref{sonsuz}.
 \hfill $\Box$ \\

\subsection{Existence results for real contact solid tori}

\begin{prop}\label{eksibir}
Up to equivariant isotopy through real tight structures relative boundary 
there is a unique $c_3$-real tight contact structure
on $S^1\times D^2$ with  convex boundary with $s=\frac{1}{2k+1}$ $(k\in \Z)$, and $\#\Gamma=2$.
\end{prop}

\noindent {\it Proof:}
The proof is the $c_3$-invariant version of the proof  of \cite[Proposition~4.3]{ho}.
Here, we outline that proof and point out required modifications.
First let us note that by sufficiently many equivariant meridional Dehn twists we can make the slope $-1$.
Hence it is enough to consider just the case $s=-1$.

By the real version of Giroux's Flexibility Theorem (Proposition~\ref{GirFlex}),
we suppose that the Legendrian rulings
have slope 0 so that each  bounds a meridional disk in $S^1\times D^2$.
Let $L$ be a Legendrian ruling and $D$ be the disk it bounds.
By \cite[Proposition~3.1]{ho}, we can perturb $D$ relative $L$ to make it convex.
Since the slope is $-1$, we have tw$(L, D)=-1$ and so the dividing set
on $D$ must be connected. It is a path with end points on the boundary.

Now,  we consider $S^1\times D^2\setminus (D \cup T^2 \cup c_{3}(D))$ which consists of two disjoint open 3-balls
that are mapped to each other via $c_3$.
Consider the closure of one of those balls. After  edge rounding and perturbing we obtain a 3-ball whose
boundary is convex and has a connected dividing set.  Up to isotopy relative boundary, there is a
unique tight contact structure on it \cite{el}. Hence define an equivariant
isotopy on $S^1\times D^2$ relative boundary
as follows: apply the isotopy on one of the 3-balls relative its boundary to obtain the model of the
unique tight structure on the ball and perform the
corresponding isotopy on the other ball antisymmetrically.
\hfill $\Box$\\

We have an interesting corollary:

\begin{coro}\label{Legkomsuluk}
Suppose $L$ is a Legendrian knot in a real tight contact manifold $(M,\xi,c)$
and in a neighborhood of $L$, $c$ is isotopic to $c_3$ through real structures. Then the twisting number of
$L$ with respect to any $c$-invariant framing is always  odd.
\end{coro}

\noindent {\it Proof:}
Take a small invariant neighborhood of $L$ with convex boundary on which $c$
acts as $c_3$. Identify the boundary with $\R^2/\Z^2$ so that the longitude determined
by the framing corresponds to $(0,1)$ while meridian corresponds to $(1,0)$.
With respect to this identification, the slope of the dividing curves on the boundary becomes
$ 1/\mbox{tw(L)}$.  It follows from Corollary~\ref{ciftyok} and Proposition~\ref{eksibir} that such
a real tight contact solid torus exists if and only if  $\mbox{tw(L)}$ is odd.
 \hfill $\Box$ \\

To cover the remaining case of  Theorem~\ref{anadolu}, we first split the solid torus
into two pieces as follows. Without loss of generality, we may assume that its core is $c_3$-invariant.
There is an invariant  Legendrian curve $L$ isotopic to the core.
In fact, we can always perturb the core to get a Legendrian curve, via creating small zigzags in the front projection.
Since each zigzag creation is local, we can create each zigzag and its symmetric counterpart 
invariantly to obtain an invariant Legendrian  $L$ at the end.

By Corollary~\ref{Legkomsuluk} the twisting number of $L$ (with respect to
the longitude $S^1\times\{(1,0)\}\subset S^1\times D^2$) is odd.
Take an invariant neighborhood $\nu(L)$ of this curve such that the boundary is invariant
and convex (Theorem~\ref{LegKomsu} and Proposition~\ref{ahenkcnvx}).
Consider the decomposition of the solid torus into two $c_3$-invariant pieces: $\nu(L)$ and the closure of
the complement of $\nu(L)$, which is the thick torus $T^2\times I$,  ($I=[0,1]$), whose boundary is the
disjoint union of two convex tori with slopes $ \pm\frac{1}{2k+1}, k\in \Z^+$  and $s$ 
respectively for inner and outer boundaries. 
Indeed, as discussed in Remark~\ref{equivdt} we can assume that  $s=-\frac{p}{q}\leq -\frac{1}{2}$ with  $p,q\in \Z^+$, $(p,q)=1$, in which case the slope of the dividing curve on $\partial \nu(L)$ is necessarily $-\frac{1}{2k+1}$. (This can be shown by an argument similar to that presented in the proof of  Lemma~\ref{kilcik} below.)

By Proposition~\ref{eksibir} we know that there is a unique real tight structure on $\nu(L)$.
In the next section our aim is to figure out the possible real tight contact structures on  
the thick torus with convex boundary with slopes $-\frac{p}{q}$ and $-\frac{1}{2k+1}$ .
The tool we use is the method of layering thick tori by means of bypasses introduced in \cite{ho}.

\subsection{Factoring $c_3$-real tight solid tori}


Now consider the real thick torus $(T^2\times I, c_3)$  with $T^2$ identified by $\R^2/ \Z^2$ as before and
suppose that we are given a real tight contact structure $\xi$ on $T^2\times I$ making the boundary convex.
We denote by $s_k$ ($k\in\{0,1\}$)  the slope of the dividing curve $\Gamma_k$ on  $T_k=T^2\times\{k\}$.
We also denote by $T(s_0,s_1)$ a $c_3$-real tight thick torus 
that has a convex boundary $T_0 \cup T_1$
with slopes $s_0$ and $s_1$ respectively and with dividing sets $\Gamma_0$ and $\Gamma_1$,
$\#\Gamma_0=2$, $\#\Gamma_1=2$. We will show that if $T^2\times I$ lies in a solid torus
then it is {\em minimally twisting}  and it can be decomposed into
{\em symmetric basic slices} and {\em symmetric double slices}
(see \cite{ho} for the definitions of minimally twisting, basic slice and bypass).
A double slice is essentially two basic slices glued together:
\begin{defi}
Given the real tight solid torus $T(s_{0},s_{1})$ as above, it is called a {\em symmetric double slice}  if on the Farey tessellation of the unit disk one can connect $s_0$ and $s_1$ in such a way that  the minimal number of edges of the paths  in the counterclockwise direction from $s_{1}$  to $s_{0}$ is two. 
\end{defi}

Let us note the following fact from \cite{ho} for the sake of completeness and for reference in the sequel.
Its proof does not require the real setting.

\begin{lemma}\label{kilcik}
Let $S^1\times D^2$ be a tight contact solid torus whose boundary is convex with slope
$s < 0$ and $\#\Gamma=2$.
Denote by  $L$  a Legendrian curve which is isotopic to the core and
by $s_{L}$ the slope of the convex boundary of the contact neighborhood
of $L$ in $S^1\times D^2$. Then $s_{L}$ lies in $[s, 0)$.
\end{lemma}

\noindent {\it Proof:}  The idea of the proof is similar to that of Proposition~4.15 in \cite{ho}. 
Assume the contrary: there is a Legendrian curve such that  $s_{L}$ is outside the interval  $[s, 0)$.   
As it is discussed in \cite{ho} by successive bypasses inwards,
we can peel off the solid torus to get solid tori with convex boundary $T_{i}$ of slope $s_{i}$.
Let  $s_{i}$ denote the slope of the convex boundary of the solid torus after $i^{th}$ bypass. Then on the Farey tessellation of the unit disk there is an edge from $s_{i}$  to $s_{i+1}$ in the counterclockwise direction 
\cite[Lemma 3.15]{ho}.  Thus a counterclockwise path from $s$ to $s_{L}$  on the Farey tessellation, obtained as a result of some number of bypasses,
necessarily passes through 0 because $s_{L}$ is outside $ [s, 0)$.  This gives a contradiction since there is no tight solid torus with slope 0.
\hfill  $\Box$ \\

Let us reconsider the real tight thick torus
$T(-\frac{1}{2k+1}, -\frac{p}{q})\subset S^1 \times D^2$.
It follows from the proof above that $-\frac{1}{2k+1} \geq -\frac{p}{q}$;
in fact otherwise there would exist a tight solid torus with slope 0.
It follows similarly that $T(-\frac{1}{2k+1}, -\frac{p}{q})$ is minimally twisting. 
By Proposition~\ref{GirFlex} we can assume that the invariant Legendrian rulings on the boundary 
tori have slope $-\frac{1}{2j}$ for some $j\in\Z^+$ with $2j>2k+1$. Consider an 
annulus $A$ bounded by a pair of Legendrian rulings 
on $T_0$ and $T_1$. By \cite[Proposition~3.1]{ho},  $A$ can be made convex relative its boundary.  
Set $S_i=\partial A \cap T_{i}$, $i=0,1$. If the difference between the intersection numbers of $S_1$ and $S_0$ with a component of 
the dividing curve of $T_{1}$   is greater than or equal to 3, namely if 
\begin{equation}
\left((q,-p)-(2k+1,-1)\right)\cdot(2j,-1)=2jp-2j+2k+1-q \geq 3,
\label{sacma} 
\end{equation} 
then there is at least a pair of bypasses that can be performed on $A$ consecutively. 
Furthermore, after performing
the first, the second bypass disk can be moved to the symmetric position of the first bypass
region so that a pair of $c_3$-symmetric bypass operations can be performed one after the other.
Given $p\neq 1$ and $q$, the inequality at (\ref{sacma}) can always be satisfied for 
a sufficiently large $j$. Note that if $p=1$, we must have $2k+1-q\geq 3$.  

For $2j>2k+1$, i.e. $-\frac{1}{2j}$ outside $[-\frac{p}{q}, -\frac{1}{2k+1}]$, 
a new pair of bypasses can always be performed 
as long as  (\ref{sacma}) is satisfied for an appropriately chosen $j$.
Note that when equality is attained then  we have a symmetric  double slice.
Hence we have proven

\begin{prop}\label{kabuklar}
Suppose $T(s_{0}, s_{1})$ $(s_0\neq s_1)$ is a real tight thick torus with $\xi$ minimally twisting
and $\#\Gamma_0=2$, $\#\Gamma_1=2$. Then $T(s_{0},s_{1})$ assumes a factoring into a number of 
symmetric double slices and possibly
a symmetric basic slice. If the minimal number of edges on the paths  from $s_{1}$ to $s_{0}$
in the counterclockwise direction on the Farey tessellation of the unit disk is odd,
there is a single symmetric basic slice in the factoring.  If that number is even, there is no
basic slice.

In particular, if $s_{1}\in\Z^-$ and  $s_0=-1$ then the last slice is either the thick torus $T(-1,-3)$ or  $T(-1,-2)$.

If $s_{1}\in\Z^-$ and $s_0=-\frac{1}{2k+1}$ $(k\in\Z^+)$ then the last slice is either $T(-\frac{1}{2k+1},-\frac{1}{2k})$
or  $T(-\frac{1}{2k+1},-\frac{1}{2k-1})$.
\end{prop}

\subsection{Some nonexistence results for $c_3$-real basic and double slices.}
Proposition~\ref{kabuklar}   suggests the investigation of the real tight structures on 
symmetric basic and double slices. This is what we are going to do now.

The peeling-off process above shows that
there is a single  basic slice occurring if the total number of bypasses performed is odd.
This basic slice has to be $T(-\frac{1}{2k+1},-\frac{1}{2k})$. However
there can be no such thick torus in the initially given solid torus since there is no solid torus 
with slope $-\frac{1}{2k}$ (Corollary~\ref{ciftyok}).

For double slices, we observe that if the total number of bypasses performed during the peeling-off process above
is even,
then the last thick torus is the double slice $T(-\frac{1}{2k+1}, -\frac{1}{2k-1})$ or
$T(-1,-3)$. For $T(-1,-3)$, we have the following

\begin{prop}\label{etlisimit3}
There is no $c_3$-real tight $T(-1,-3)$.
\end{prop}
\noindent {\it Proof:}
We assume that there exists a real contact thick torus with the required properties
and conclude that in each case possible, the contact structure must be overtwisted. 
Suppose the slope of Legendrian rulings is set to be $0$ on both $T_0$ and $T_1$.
Let $A$ denote the convex annulus obtained by a small perturbation relative boundary 
of an annulus bounded by Legendrian rulings on $T_0$ and $T_1$. 
Set $S_1=\partial A \cap T_1$ and $S_2=\partial A \cap T_2$. 
We observe that $\mbox{tw}(S_1,T_1)=-3$ and $\mbox{tw}(S_0,T_0)= -1$
Therefore on $A$ there are at least two
arcs of the dividing set which start and end at $S_1$ \cite[Proposition~3.17]{ho}.

The real structure $c_3$ gives a disjoint copy $\bar{A}=c_{3}(A)$ of $A$
which has the same characteristic foliation (and the same dividing set) as $A$ except
the signs of the singularities and the orientation of the foliation reversed.
The complement of $\bar{A}\cup A$ consists of two connected components,
each of which  has as closure a topological solid torus. These solid tori are symmetric with respect to $c_3$.
Let $S$ denote one of those solid tori. Then we consider the contact solid torus $\hat{S}$ obtained by
rounding the edges of $S$ \cite[Lemma~3.11]{ho}.
Our aim is to investigate the contact structure on $\hat{S}$ by means of the dividing set of its boundary.
While rounding the edges the dividing set is constructed by going towards right at each surface change and following the dividing curve first encountered \cite{ho}. 

\begin{figure}[h]
\begin{center}
\resizebox{10cm}{!}
{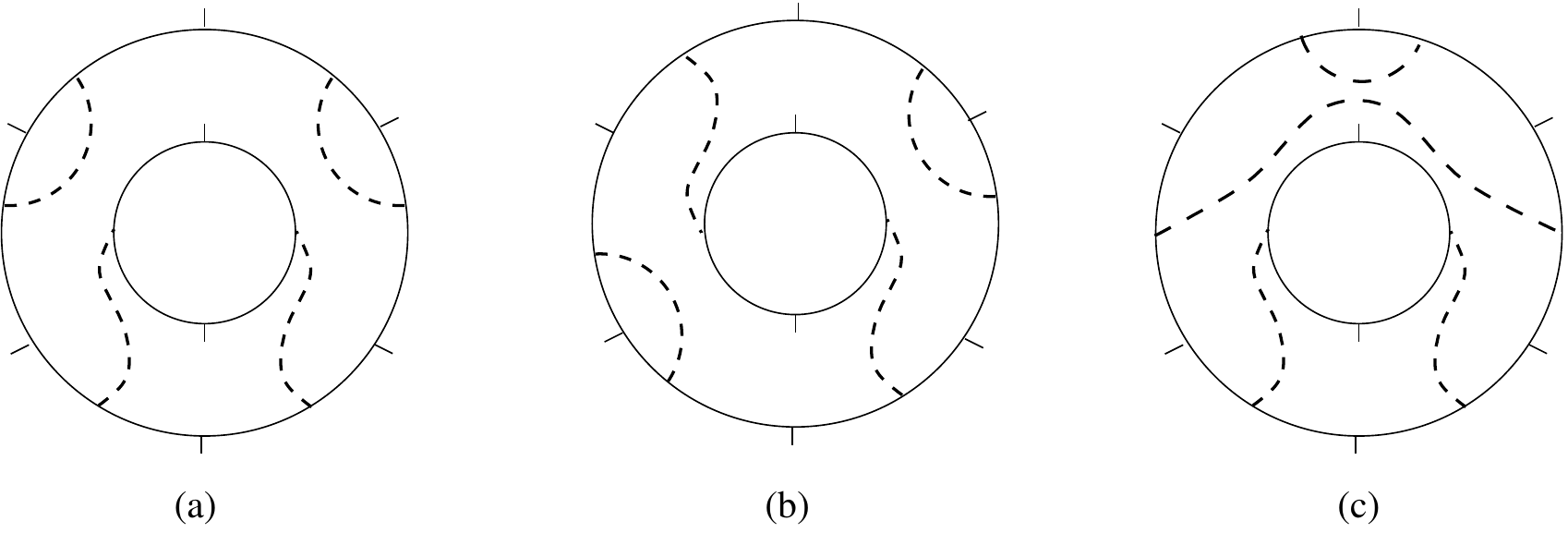}
\caption{Possible dividing sets on $A$ in $T(-1, -3)$}
\label{ciplaklarplaji3}
\end{center}
\end{figure}

We have several possible configurations of the dividing set on the annulus $A$.
First let us note that  the dividing set on $A$ can not have a closed 
component which bounds a disk on $A$ since in that case  the contact structure would be overtwisted
\cite{gi2}.
Note  also that by an argument similar to  Step 1 of the proof of \cite[Proposition~4.7]{ho}  one can show that 
every dividing arc starting from the inner boundary of $A$ must end at
the outer boundary of $A$. Furthermore, those arcs can be put in a {\em standard}  position
in which the number of their rotations around the inner boundary can be made 0
(this is called {\em holonomy}  in \cite{ho}; there the interior of the chosen annulus is pushed
one turn around the torus to change the holonomy by $\pm 1$.)
We do not present how to attain 0 holonomy in or setting since  
the discussion below is insensitive to the holonomy.
Those said, there are three possible configurations for the dividing curve on $A$
(up to  holonomy), depicted in Figure~\ref{ciplaklarplaji3}.

\begin{figure}[h]
\begin{center}
\resizebox{5.3cm}{!}
{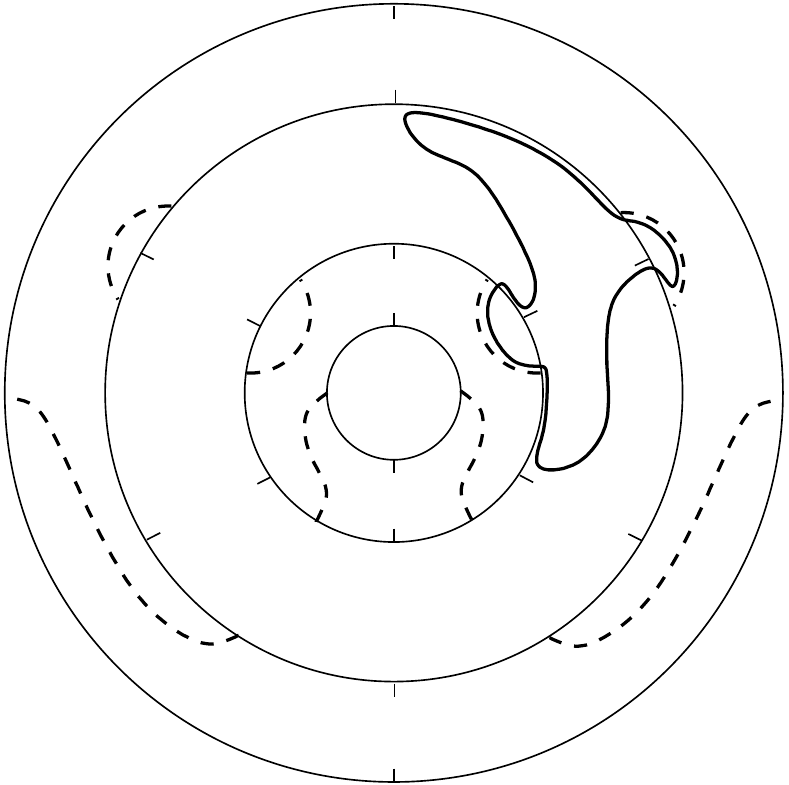}
\caption{The dividing set on $\hat{S}$ for the configuration in the case (a).}
\label{hicyok}
\end{center}
\end{figure}

\begin{figure}[h]
\begin{center}
\resizebox{5.3cm}{!}
{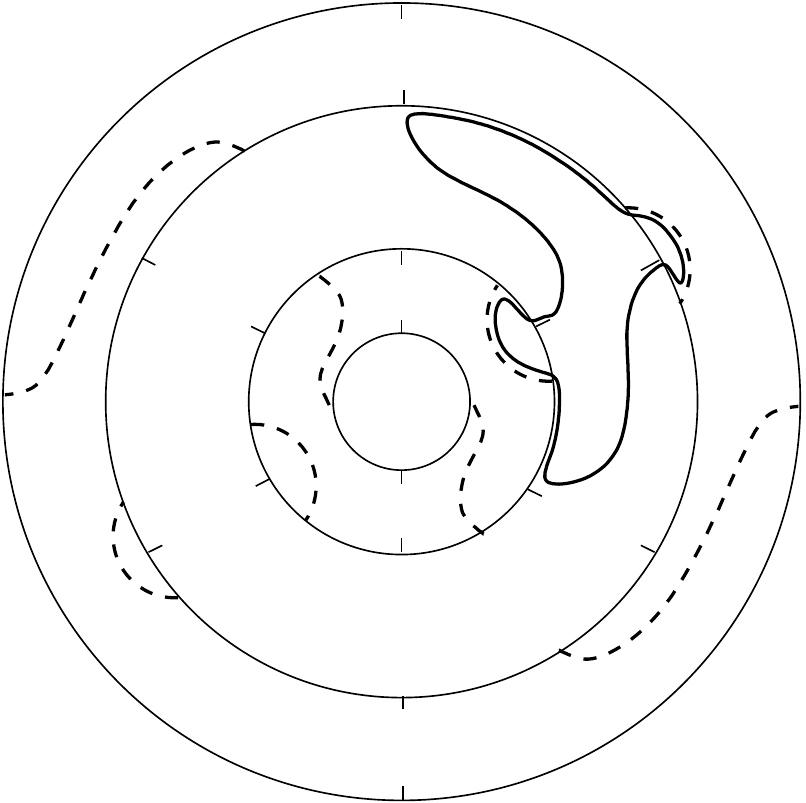}
\caption{The dividing set on $\hat{S}$ for the configuration in the case (b).}
\label{sifirsimit}
\end{center}
\end{figure}

With those configurations, Figure~\ref{hicyok}  proves that there is no real tight  
$T(-1,-3)$ leading to the configuration in Figure~\ref{ciplaklarplaji3}(a). 
Indeed,  we observe that the contact structure around $\hat{S}$ must be overtwisted
since the dividing curve in Figure~\ref{hicyok} obtained after edge-rounding bounds a disk on 
the outer boundary  torus $T_1$. 
In Figure~\ref{hicyok} the innermost and outermost annuli depict $A$ and $\bar{A}$ respectively. 
The annulus in between shows half of  $T_1$. Slope on $T_1$ is $-3$, therefore going from $A$ to $\bar{A}$
the dividing curve rotates $-\frac{1}{6}$ of full turn.
The  annulus corresponding to half of $T_0$ is on the reverse of the figure, which we do not see. 
The torus in the figure thus obtained is oriented in accordance with the orientation of $T_1$.

Simlarly, Figure~\ref{sifirsimit}  proves that there is no real tight  
$T(-1,-3)$ leading to the configuration in Figure~\ref{ciplaklarplaji3}(b).
Case (c) is similar and we do not draw a figure for it. \hfill  $\Box$

\begin{figure}[h]
\begin{center}
\resizebox{8.2cm}{!}
{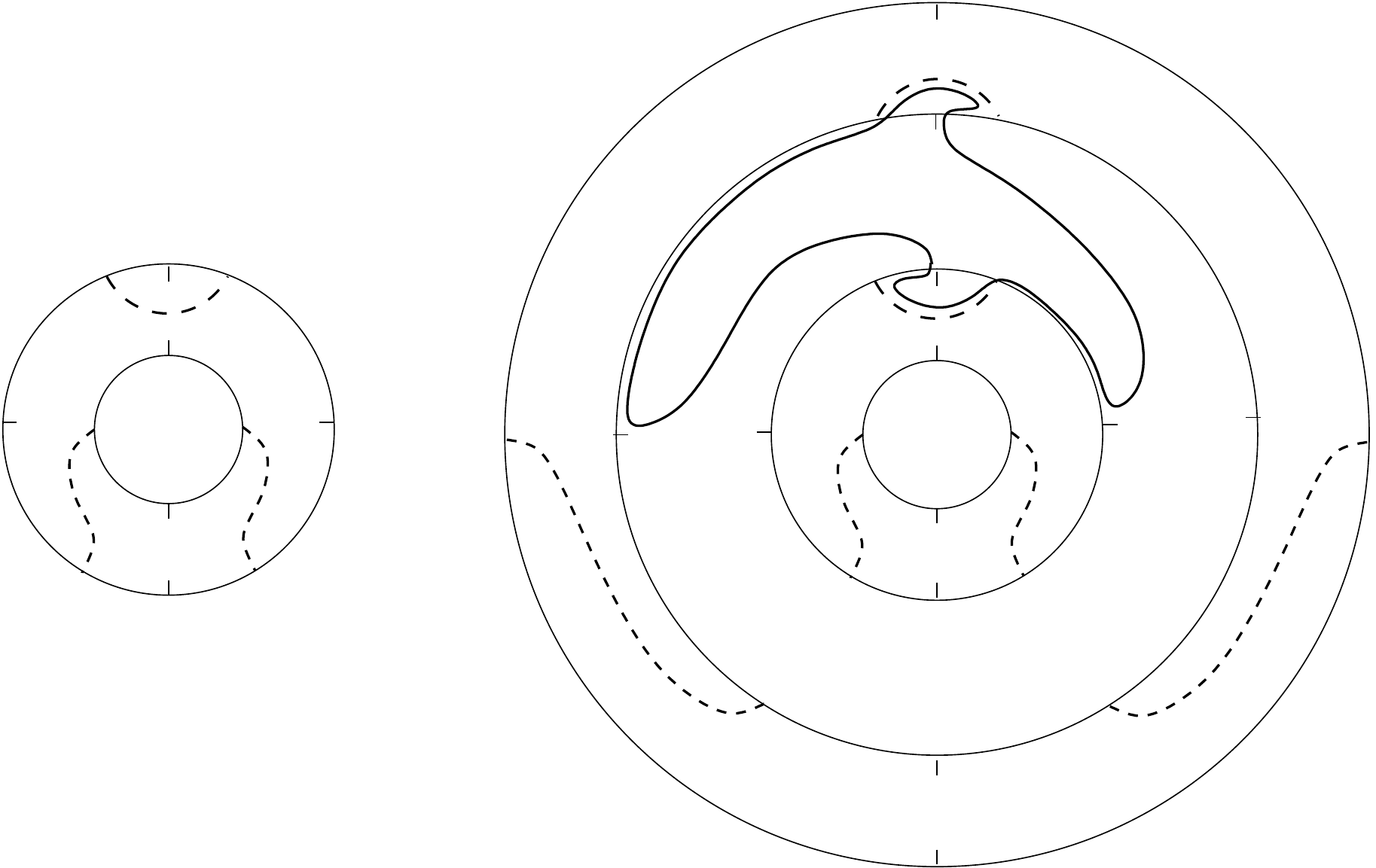}
\caption{The single possible configuration (left) of the dividing curve on $A$ in $T(-1, -2)$ 
begets an overtwisted disk (right).}
\label{otsimit}
\end{center}
\end{figure}

In a similar way one can prove the following
\begin{prop}\label{etlisimit2}
There is no $c_3$-real tight $T(-1,-2)$.
\end{prop}
\noindent {\it Proof:}
The proof is exactly as the previous. The only 
configuration on $A$ up to holonomy is as in Figure~\ref{otsimit}, which leads an overtwisted disk as shown.
\hfill  $\Box$\\

\noindent {\em Proof of Theorem~\ref{anadolu}:}

Given a $c_3$-real solid torus with slope $-\frac{p}{q}$,
the case $p=1$ is covered by Corollary~\ref{ciftyok} and Proposition~\ref{eksibir}.
If $q=1$ and $p>1$ then  we split the solid torus into two pieces as in the discussion
following Corollary~\ref{Legkomsuluk}.
An invariant isotopic  Legendrian copy  $L$ of the core has a standard symmetric tight neighborhood
with the boundary invariant and convex.
We decompose the solid torus into two invariant pieces: $\nu(L)$ and  $T(-\frac{1}{2k+1},-\frac{p}{1})$.
Following Lemma~\ref{kilcik}, the latter was observed to be minimally twisting.
After  $c_3$-invariantly factoring the thick torus, among the slices we obtain either 
a final slice $T(-\frac{1}{2k+1},-\frac{1}{2k})$ or $T(-1,-2)$
or an intermediate slice $T(-1,-3)$  (the first two and the last cases are mutually exclusive).
None of these cases is possible thanks to the first paragraph of this section,  Proposition~\ref{etlisimit2} and Proposition~\ref{etlisimit3} respectively.
\hfill  $\Box$

\end{document}